\documentclass[reviewcopy]{elsart}
\usepackage{amssymb,amsmath,url,color}
\usepackage[mathscr]{eucal}
\def\mhat#1#2#3{T_{#1,#2}(#3)}
\def\re{\mathop{\rm Re}\nolimits}

  \def\gg#1{\Gamma(#1)}

  \def\nn{\nonumber}
  \def\nr{\nn\\}
  \def \ds{\displaystyle}
  \def\res#1#2{\hbox{\rm Res}_{#1}\,#2}
  \def\wwhit{W_{n,a}(z)}
  
  \def\mhat#1#2#3{T_{#1,#2}(#3)}
  \def\gam#1{\Gamma(\displaystyle{#1})}
  
  \def\R{\mathbb{R}}
   \def\C{\mathbb{C}}  
   \def\Z{\mathbb{ Z}}  
   
\begin{document}
\begin{frontmatter}
\title{Recurrence relations for Mellin transforms\\   of $GL(n,\R)$ Whittaker functions} 
\author{Eric Stade}
\address{Department of Mathematics, University of Colorado, Boulder, CO 80309, USA }
 \ead{stade@colorado.edu} 
\author{Tien Trinh} 

\address{Department of Mathematics, Hanoi National University of Education, 136 Xuan Thuy St., Hanoi, Vietnam }

\ead{tientd@hnue.edu.vn}

\begin{abstract}Using a recursive formula for the Mellin transform $T_{n,a}(s)$ of a spherical, principal series $GL(n,\R)$ Whittaker function, we develop an explicit recurrence relation for this Mellin transform. This relation, for any $n\ge2$, expresses $T_{n,a}(s)$ in terms of a number of ``shifted''  transforms $T_{n,a}(s+\Sigma)$, with each coordinate of $\Sigma$ being a non-negative integer.

 We then focus on the case $n=4$. In this case, we use the relation referenced above to derive further relations, each of which involves ``strictly positive shifts'' in one of the coordinates of $s$.   More specifically:  each of our new relations expresses  $T_{4,a}(s)$ in terms of  $T_{4,a}(s+\Sigma)$ and $T_{4,a}(s+\Omega)$, where for some $1\le k\le 3$, the $k$th coordinates of both $\Sigma$ and $\Omega$ are strictly positive.  
 
 Next, we deduce a recurrence relation for $T_{4,a}(s)$ involving strictly positive shifts in all three $s_k$'s at once.  (That is, the condition ``for some $1\le k\le 3$'' above becomes ``for all $1\le k\le 3$.'')
 
 These additional relations on $GL(4,\R)$ may be applied to the explicit understanding of certain poles and residues of $T_{4,a}(s)$. This residue information is, as we describe below, in turn relevant to recent results concerning orthogonality of Fourier coefficients of $SL(4,\Z)$ Maass forms, and the $GL(4)$ Kuznetsov formula.
  \end{abstract}

  \end{frontmatter}

  \section{Introduction; definitions and notation}
Friedberg and Goldfeld, in \cite{FG},  demonstrated the existence of recurrence relations for the Mellin transform $T_{n,a}(s)$ of a spherical, principal series $GL(n,\R)$ Whittaker function.  Additionally, they described a method for producing certain relations of this kind explicitly, given explicit expressions for the $GL(n,\R)$-invariant differential operators on the generalized upper half-plane $\mathscr{H}_n=GL(n,\R)/(O(n,\R)\times \R^*)$.

In the present work, we develop such relations by a different method, which does not require explicit knowledge of the differential operators on $\mathscr{H}_n$, and which provides relations more readily, and in more manageable forms. Our approach uses a recursive expression, developed in \cite{IshiiStade1}, for $T_{n,a}(s)$ in terms of $T_{n-1,a'}(s')$ (for certain $a'$ and $s'$ depending on $a$ and $s$).

We then focus on the case $n=4$.  We combine our above recurrence relation, in this case, with translates of this relation under a certain action of the Weyl group for $GL(4,\R)$ on the variable $a$.  We thereby obtain new relations that are particularly useful for studying meromorphic continuation, poles, and residues of $T_{4,a}(s)$.

Such information is, in turn, relevant to recent work \cite{GW}  concerning orthogonality of Fourier coefficients of Maass forms for $GL(4)$. This work extends results previously developed in the context of $GL(3)$, cf. \cite{GK} (particularly Theorem 1.3 there).  In the $GL(4)$ context, the analyses entail a test function $p_{T,R}(y_1,y_2,y_2)$ depending on parameters $T$ and $R$.    In \cite{GW},   $p_{T,R}$ is defined  in terms of its explicitly specified Lebedev-Whittaker transform $p_{T,R}^\sharp$ (cf. \cite{GK2}, \cite{Wallach2}). Because the former equals the integral of the latter against a spherical, principal series $GL(4,\R)$ Whittaker function, integrals -- specifically, Mellin transforms -- of these Whittaker functions arise.  To understand $p_{T,R}$ (especially its growth properties), then, one must have information concerning various properties of $T_{4,a}(s)$.  We will investigate these properties  in this paper.

To establish the framework for our results, and to provide detail concerning the ideas discussed above, let us first
discuss harmonic analysis on
$GL(n,\R)$. To this end, we let
$X_n\subset
GL(n,\R)$ be the group of upper triangular matrices with
diagonal elements equal to one, and let
$Y_n\subset
GL(n,\R)$ be the group of diagonal matrices $y$ of the form
\[y=\hbox{diag}(y_1y_2\cdots y_{n-1},y_2y_3\cdots
y_{n-1},\ldots ,y_{n-1},1),\]where $ y_j\in\R^+$ for $1\le
j\le n-1.$ The Iwasawa decomposition for $GL(n,\R)$ identifies the
``generalized upper half-plane''    \[
\mathscr{H}^n=GL(n,\R)/(O(n,\R)\cdot{\R^\ast})   \] with the set 
    \[\lbrace z= x y\,|\,
x\in X_n, y\in Y_n\rbrace.   \]   

We wish to consider   certain (non-zero) eigenfunctions of the algebra 
$D$ of $GL(n,\R)$-invariant differential operators on $\mathscr{H}^n $.
To do so, it is convenient to start with the simplest such
eigenfunction, namely,   the {\it power
function} 
\begin{equation}\label{pfun} \prod_{j=1}^{n-1} y_j^{ j(n-j)/2}\, 
\prod_{k=1}^{j}{y_j}^{2a_k},\end{equation} where $a_k\in\C$ for
$1\le j\le n-1$.  Let us put
$a_n=-a_1-a_2-\cdots-a_{n-1}$ and
$a=(a_1,a_2,\ldots,a_n)$, and denote the corresponding power function 
\eqref{pfun} by $H_a(z)\equiv H_a(y)$. It is then shown \cite{Sel} that
$H_a$ is indeed an eigenfunction of $D$, and that its
eigenvalues $\lambda_a(d)$, defined by
    \[dH_a=\lambda_a(d)H_a\quad (d\in D),   \]  are invariant under any
permutation of the $a_j$'s. 

$GL(n,\R)$ Whittaker functions may now be defined in terms of the
above power function $H_a$, and the character
    \[\Theta(x)=\hbox{exp}(2\pi i(x_{1,2}+x_{2,3}+\cdots
+x_{n-1,n}))   \] of $X_n$. We have:
\begin{defn} \label{whitdef}A $GL(n,\R)$  Whittaker function of type
$a$ is a function $f_a$, smooth on
$\mathcal{ H}^n$
and meromorphic in $a_1,a_2,\ldots,a_{n-1}$, such that:

{(a)}
$df_a=\lambda_a(d)f_a$ for all $d  \in D$;

{(b)} $f_a(x'
z) =\Theta(x')f_a(z)$ for all $x' 
\in X_n,  z \in \mathcal{ H}^n$.
\end{defn}

The systematic study of Whittaker functions was
initiated, in the more general context of Chevalley groups over
local fields, by Jacquet \cite{Jacquet}.  We will restrict our attention here to the context of $GL(n,\R)$.

The space
$V_{n,a}$ of
$GL(n,\R)$ Whittaker functions of type $a$ has dimension
$n!$: this  follows from unpublished work of Casselman and
Zuckerman, and independent work of  Kostant
\cite{Kostant}.  In this paper, we will be concerned with the  spherical
principal series  Whittaker function
$\wwhit\in V_{n,a}$, which is (up to scalars) the unique element of
$V_{n,a}$ that decays rapidly as $y_j\rightarrow\infty$
for
$1\le j\le n-1$. 
This uniqueness property of $\wwhit$ follows
from multiplicity-one theorems of Shalika  \cite{Shalika} and Wallach
\cite{Wallach}.

The Whittaker function $\wwhit$ is central to the harmonic analysis of 
automorphic forms on $GL(n,\R)$.  Specifically,
suppose  $\varphi\colon \mathcal{ H}^n\rightarrow \C$ is a {\it cusp form} for $GL(n,\Z)$ -- in other words,
$\varphi(\gamma z)=\varphi(z)$ for all $\gamma\in GL(n,\Z)$ and
$z\in\mathcal{ H}^n$;  $\varphi$ decays rapidly in each $y_j$,
as $y_j\rightarrow\infty$; and $\varphi$ is an eigenfunction of
$D$ with $d\varphi=\lambda_a(d)\varphi$ for all $d$.  Then $\varphi$ has a ``Fourier-Whittaker  expansion'' (cf.
\cite{Piatetski},  \cite{Shalika}), all of whose terms are expressible in terms of
$\wwhit$.   This follows from the multiplicity-one theorems mentioned above.  If $\varphi$ is instead a {\it Maass form}, which is similar to a cusp form but may  satisfy somewhat less stringent growth conditions, 
it is still the case that ``most'' 
of the Fourier coefficients of $\varphi$ {are} expressible in
terms of $\wwhit$.
 
The (normalized)  {\it Mellin transform} $\mhat{n}{a}{s}$ of $W_{n,a}(z)$, given by
\begin{equation}\label{melltran}\mhat{n}{a}{s}=2^{n-1}\int_0^\infty\ldots\int_0^\infty
W_{n, a}(y_1,y_2,\ldots,y_{n-1})\biggl\{\frac{1}{y^{\rho_n}} \prod_{j=1}^{n-1}(\pi
y_j)^{2s_j} \,\frac{dy_j}{
y_j}\biggr\}\end{equation}for $s=(s_1,s_2,\ldots,s_{n-1})\in {
\C}^{n-1}$ and   \begin{align} \label{yrhon} y^{\rho_n} = \prod_{j=1}^{n-1} y_j^{j(n-j)/2},\end{align}    plays a significant role in the theory of automorphic
forms, and especially in the study of automorphic $L$-functions.   (See, for example, \cite{Bump87}, \cite{JPS}, \cite{JS}, \cite{Moreno}, \cite{Shimura},  \cite{Stade01}, \cite{Stade02}.)  We remark that $\wwhit$, properly normalized, is {\it invariant} under
permutations of the $a_j$'s, and consequently so is
$\mhat{n}{a}{s}$. The standard normalization of
$\wwhit$, which is the one originally given by Jacquet \cite{Jacquet},  and to
which we adhere  throughout this paper, ensures this
invariance.
 
By Mellin inversion, we have
\begin{align} \label{Wdef}
  W_{n,a}(y) & = W_{n,a}(y_1, \cdots, y_{n-1} ) = 
\frac{y^{\rho_n} }{ (2\pi i)^{n-1}}
\int_{s_1,\cdots,s_{n-1}}\mhat{n}{{a}}{s}
\prod_{j=1}^{n-1}(\pi y_j)^{-2s_j} \,ds_j, 
\end{align}   with 
the path of integration in each
$s_j$ being a vertical line in the complex plane, of sufficiently large real part  to keep the poles of
$\mhat{n}{{a}}{s}$ on its left.

We conclude this section with  some basic facts --  to be of importance in what follows --  concerning gamma
and  Bessel functions. Proofs of all of
these facts may be found, for example, in \cite{WhWa}.

For the gamma
function,  defined by
\begin{align}  \gam{s}=\int_0^\infty
\hbox{e}^{-t} \,t^s\,\frac{dt}{t} \nonumber \end{align}   for $\re{s}>0$, and having  analytic continuation to $\C\backslash\{0,-1,-2,\ldots\}$, we have the
residue formula\begin{align}  \hbox{\rm Res}(\gam{s},s=-n)=\frac{(-1)^n}{
n!}\quad(n=0,1,2,\ldots) ,  \nonumber  \end{align}    and the translation and reflection  formulas 
\begin{align}  \label{gammaids}\Gamma(s+1)=s\Gamma(s),\quad \Gamma(s)\Gamma(1-s)=\frac{\pi}{\sin\pi s}.  \end{align}    

Next, the {modified} Bessel function  $K_\nu(z)$ of the second kind is defined as follows: 
\begin{align}  \label{Kbess2}K_\nu(z)=\frac{1}{4}\cdot \frac{1}{2\pi i} \int_s\Gamma\biggl(\frac{s+\nu}{2}\biggr)\Gamma\biggl(\frac{s-\nu}{2}\biggr)\biggl(\frac{z}{2}\biggr)^{-s}\,ds,  \end{align}   the path of integration being a vertical line to the right of any poles of the integrand; or, equivalently\begin{align}  \label{Kbess3}K_\nu(z)=\frac{1}{2} \int_0^\infty e^{-y(t+t^{-1})/2} \,t^\nu \frac{dt}{t}\qquad(\re{y}>0).  \end{align}   

As is well known (see, for example, \cite{Bump}), the spherical principal series Whittaker function on $GL(2,\R)$ is closely related to the Bessel  function $K_\nu(z)$    -- specifically, we have
\begin{equation}W_{2,a}(y)=2\sqrt{y}K_{2a}(2\pi y).\label{W2}\end{equation} (Here, by a slight abuse of notation, we identify $(a,-a)\in  \C^2$ with the complex number $a$, and identify the matrix $y=\hbox{diag}(y,1)\in Y_n$   with the positive number $y$.)    
    \section{A recurrence formula for general \boldmath $n$}
We begin with the following recursive formula, derived in \cite{IshiiStade1}, for the Mellin transform $T_{n,a}(s)$: \begin{align} &T_{n,a}(s)=\frac{\Gamma(s_1+a_1)\Gamma(s_{n-1}-a_1)}{(2\pi i)^{n-2}}\nn\\&\times\int_{z}\Bigg(\prod_{j=1}^{n-2}\Gamma\Bigl(z_j+s_j-\frac{ja_1}{n-1}\Bigr)\Gamma\Bigl(z_j+s_{j+1}+\frac{(n-1-j)a_1}{n-1}\Bigr)\Bigg)\, T_{n-1,b}(-z)\,dz\nn\\&=\frac{\Gamma(s_1+a_1)\Gamma(s_{n-1}-a_1)}{(2\pi i)^{n-2}}\nn\\&\times\int_{z}\Bigg(\prod_{j=1}^{n-2}\Gamma(z_j+s_j)\Gamma(z_j+s_{j+1}+a_1)\Bigg)\, T_{n-1,b}\biggl(\Big(-z_j-\frac{ja_1}{n-1}\Big)_{1\le j\le n-2}\biggr)\,dz,\label{Tnas}\end{align}  
where $a=(a_1, a_2,\ldots, a_n)$  and $s=(s_1,s_2,\ldots, s_{n-1})$, as in the previous section; also,$$\displaylines{ z=(z_1,z_2,\ldots, z_{n-2}),\quad dz=dz_1\,dz_2\cdots dz_{n-2} ,\cr b=(b_1,b_2,\ldots, b_{n-1})=\biggl(a_2+\frac{a_1}{n-1}, a_3+\frac{a_1}{n-1},\ldots, a_{n}+\frac{a_1}{n-1}\biggr).}$$(The second equality in \eqref{Tnas} follows from the change of variable $z_j\to z_j+ja_1/(n-1)$ applied to the first.)
It is shown in \cite{ST1} that the above formula defines $T_{n,a}(s)$ as an analytic function of $s$, in a product of half-planes where the real parts of the $s_k$'s are sufficiently large with respect to the real parts of the $a_j$'s.

\begin{rem}{\rm The right-hand side of \eqref{Tnas} makes sense for $n\ge2$, provided we   agree that the ``zero-fold'' integral that appears there, in this case, simply equals one.  With these conventions in force, \eqref{Tnas} implies that\begin{equation} T_{2,a}(s)=\gg{s_1+a_1}\gg{s_1-a_1},\label{gl2}\end{equation}a result that also follows   from \eqref{Kbess2},   \eqref{W2}, and the transform pair \eqref{melltran}, \eqref{Wdef}.}\end{rem}

In this section, we will use 
\eqref{Tnas} to deduce a different sort of recurrence relation for $T_{n,a}(s)$ -- one that expresses this Mellin transform in terms of ``translates'' $T_{n,a}(s+\Sigma)$, for certain integer $(n-1)$-tuples $\Sigma$.  The main idea behind our derivation of this relation is as follows.  Because of the functional equation $\gg{s+1}=s\gg{s}$, replacing $s$ by $s+\Sigma$ in  \eqref{Tnas} will result in some ``extra'' factors in the integrand on the right-hand side.  If some linear combination of these factors, with coefficients independent of the $z_j$'s, yields a nonzero expression that is also independent of the $z_j$'s, then we have a representation of $T_{n,a}(s)$ in terms of the corresponding translates $T_{n,a}(s+\Sigma)$.

To this end, we have the following combinatorial lemma.
\begin{lem}
 Let $n$ be an integer with $n\ge 2$.  Let\[U_{n-1}=\biggl\{\mu=(\mu_1,\cdots, \mu_{n-1})\in \{0,1\} ^{n-1} \hbox{ \rm and } \sum_{i=1}^{n-2}\mu_i\mu_{i+1}=0\biggr\}\] be the set of binary sequences of length $n-1$ that do not contain two adjacent 1's, and define 
\[I_{\mu}=\{i\in \{1,2,\cdots,n-1\}\colon\mu_i=1\}\hspace{0.5cm}\qquad(\mu\in U_{n-1}).\]   
Let  $z_0=s_0=s_n=0$ and $z_{n-1}=-a_1$, and define
\begin{equation}\alpha_k=\frac{(z_{k-1}+s_k+a_1)(z_{k}+s_k)}{(s_{k-1}-s_k-a_1)(s_{k}-s_{k+1}-a_1)}.\label{alphadef}\end{equation}Then (defining the empty product to equal $1$) we have\begin{equation}\sum_{\mu\in U_{n-1}}\prod_{k\in I_\mu}\alpha_k=0.\label{lemmaC}\end{equation}\label{thislemma}\end{lem}
\begin{pf}  We proceed by induction on $n$.  Since $U_{1}=\{(0),(1)\}$, we see that
\begin{equation*}\sum_{\mu\in U_{1}}\prod_{k\in I_\mu}\alpha_k=\prod_{k\in I_{(0)}}\alpha_k+\prod_{k\in I_{(1)}}\alpha_k=\prod_{k\in \emptyset}\alpha_k+\alpha_1=1+\frac{( s_1+a_1)(-a_1+s_1)}{( -s_1-a_1)(s_{1} -a_1)}=0, \end{equation*}so \eqref{lemmaC} is true in the base case $n=2$.
Now assume that \eqref{lemmaC} is true for a given integer $n\ge2$.  For clarity of notation, we write $\beta_k$ for the quantity defined on the right-hand side of \eqref{alphadef}, but with $n+1$ in place of $n$.  It follows readily that
\begin{equation}\beta_k=\alpha_k(1+\delta_{k,n-1}\beta_n)\quad\hbox{for } 1\le k\le n-1,\label{betadef}\end{equation}where $\delta_{k,n-1}$ denotes the usual Kronecker delta function. To complete our induction proof, we need to show that \begin{equation}\sum_{\mu\in U_{n}}\prod_{k\in I_\mu}\beta_k=0.\label{ind-step}\end{equation} 
No $\mu\in U_{n}$ can end in a pair of $1$'s, so we may write
\begin{equation}\sum_{\mu\in U_{n}}\prod_{k\in I_\mu}\beta_k=\sum_{\mu\in U_{n}\atop \mu_{n-1}=0}\prod_{k\in I_\mu}\beta_k +\sum_{\mu\in U_{n}\atop (\mu_{n-1},\mu_n)=(1,0)}\prod_{k\in I_\mu}\beta_k.\label{breakup}\end{equation}
We now investigate each of the two sums on the right-hand side of \eqref{breakup}.
 On the one hand, since $\beta_k=\alpha_k$ for $1\le k\le n-2$, the first of these two sums may be rewritten as follows:
 \begin{align} \sum_{\mu\in U_{n}\atop \mu_{n-1} =0}\prod_{k\in I_\mu}\beta_k&
=\sum_{\mu\in U_{n}\atop (\mu_{n-1},\mu_n)=(0,0)}\prod_{k\in I_\mu}\beta_k
+\sum_{\mu\in U_{n}\atop (\mu_{n-1},\mu_n)=(0,1)}\prod_{k\in I_\mu}\beta_k  
\nr&
=\sum_{\mu\in U_{n-1}\atop \mu_{n-1}=0}\prod_{k\in I_\mu}\alpha_k
+\sum_{\mu\in U_{n-1}\atop \mu_{n-1}=0}\beta_n\prod_{k\in I_\mu}\alpha_k =(1+\beta_n)\sum_{\mu\in U_{n-1}\atop \mu_{n-1}=0}\prod_{k\in I_\mu}\alpha_k\nr& =-(1+\beta_n)\sum_{\mu\in U_{n-1}\atop \mu_{n-1}=1}\prod_{k\in I_\mu}\alpha_k  \label{breakup12},\end{align}the last step by the induction hypothesis \eqref{lemmaC}.
On the other hand, the second sum in \eqref{breakup} may be rewritten using \eqref{betadef}, together with the fact that, by appending a zero to each element of $U_{n-1}$ that ends in a $1$, we obtain exactly those elements of $U_{n}$ that end in $1,0$.  Thus
\begin{align}\sum_{\mu\in U_{n}\atop (\mu_{n-1},\mu_n)=(1,0)}\prod_{k\in I_\mu}\beta_k= \sum_{\mu\in U_{n-1}\atop \mu_{n-1}=1}\prod_{k\in I_\mu}\beta_k=(1+\beta_n) \sum_{\mu\in U_{n-1}\atop \mu_{n-1}=1}\prod_{k\in I_\mu}\alpha_k.\label{breakup3}\end{align}Putting \eqref{breakup12} and  \eqref{breakup3} into  \eqref{breakup} gives 
\begin{equation*}\sum_{\mu\in U_{n}}\prod_{k\in I_\mu}\beta_k=0,\end{equation*}and we are done.
\hfill  $\Box$\end{pf} 

\begin{rem}{ \rm If we append $0,1$ to each element of $U_{n-2}$,  append $0$ to each element of $U_{n-1}$, and take the union of the results, we get exactly the set $U_n$.  Thus $|U_n|=|U_{n-1}|+|U_{n-2}|$.  From this, and the facts that $|U_1|=|\{(0),(1)\}|=2$ and $|U_2|=|\{(0,0),(0,1),(1,0)\}|=3$, we see that the number of summands in \eqref{lemmaC} (and hence in \eqref{bigthm} below) equals the $(n+1)$st Fibonacci number $F_{n+1}$.}\end{rem}

Now let \[m_k=\frac{1}{(s_{k-1}-s_k-a_1)(s_k-s_{k+1}-a_1)}\quad (k=1,2,\ldots, n-1),\]where again $s_0=s_n=0$.  The main result of this section is the following.
\begin{thm}\label{mainthm}
\begin{equation}\sum_{\mu\in U_{n-1}}\bigg(\prod_{k\in I_\mu}m_k\bigg)  \,T_{n,a}(s+\mu)=0.\label{bigthm}\end{equation}\end{thm}
\begin{pf} By \eqref{Tnas} and the functional equation $\gg{s+1}=s\gg{s}$, the left-hand side of \eqref{bigthm} equals
\begin{align*}&\frac{\Gamma(s_1+a_1)\Gamma(s_{n-1}-a_1)}{(2\pi i)^{n-2}}\int_{z}\bigg(\sum_{\mu\in U_{n-1}}\prod_{k\in I_\mu}\frac{(z_{k-1}+s_k+a_1)(z_{k}+s_k)}{(s_{k-1}-s_k-a_1)(s_{k}-s_{k+1}-a_1)}\bigg) 
\\&\times\Bigg(\prod_{j=1}^{n-2}\Gamma(z_j+s_j)\Gamma(z_j+s_{j+1}+a_1)\Bigg)\, T_{n-1,b}\biggl(\Big(-z_j-\frac{ja_1}{n-1}\Big)_{1\le j\le n-2}\biggr)\,dz,\end{align*} 
which is zero by Lemma \ref{thislemma}.
\hfill  $\Box$\end{pf} 

\begin{rem}{\rm In the case $n=2$, \eqref{bigthm} yields \begin{equation}T_{2,a}(s)=\frac{T_{2,a}(s+1)}{(s_1+a_1)(s_1-a_1)},\label{recur2}\end{equation}a result that also follows from \eqref{gl2} and the equation $\gg{s+1}=s\gg{s}$.  We may consider \eqref{recur2} to be something of a prototype for  the various recurrence relations that follow.  Compare, for example, with Propositions \ref{firstrecur}, \ref{secondrecur}, \ref{thirdrecur}, and \ref{fourthrecur} below.}\end{rem}
\section{The case $n=4$}
We now restrict our attention to the case $n=4$, and investigate some analytic properties of $T_{4,a}(s)$.  Note that Theorem \ref{mainthm} reads, in this case,
\begin{align}  &T_{4,a}(s_1,s_2,s_3)+\frac{T_{4,a}(s_1+1,s_2,s_3)}{( -s_1-a_1)(s_1-s_{2}-a_1)}+\frac{T_{4,a}(s_1,s_2+1,s_3)}{(s_{1}-s_2-a_1)(s_2-s_{3}-a_1)}
\nn\\&+\frac{T_{4,a}(s_1,s_2,s_3+1)}{(s_{2}-s_3-a_1)( s_{3}-a_1)}+\frac{T_{4,a}(s_1+1,s_2,s_3+1)}{( -s_1-a_1)(s_1-s_{2}-a_1)(s_{2}-s_3-a_1)( s_{3}-a_1)} \nn\\&=0.\label{gl4rec}\end{align}  
From this relation, we wish to derive new identities that relate our Mellin transform $\mhat{4}{a}{s}$ to shifts of that transform that are ``strictly positive'' in one or more of the variables $s_1$, $s_2$, and $s_3$.  Such identities will be useful for deriving meromorphic continuation and residue properties of our transform, cf. Section 4 below.

We begin with a pair of relations entailing strictly positive shifts in $s_1$.  Each of these relations also involves a shift in one of the other $s_k$'s -- in the first case, $k=3$, and in the second, $k=2$.

\begin{prop} \label{firstrecur}Suppose  $s_1\ne -a_k$ for $1\le k\le 4$.

\noindent {\rm (a)} We have \begin{equation}\mhat{4}{a} {s_1, s_2, 
    s_3} =\frac{ B(s_1,s_2,s_3) \mhat {4}{a} {s_1 + 1, s_2, 
    s_3}+\mhat {4}{a}{s_1 + 1, s_2, s_3 + 1}  }{\prod_{k=1}^4(s_1+a_k)},  \label{prop1} \end{equation}where
     \begin{equation}B(p,q,r)=(p+q-r)(p+r). \label{Bdef}\end{equation} 
    
\noindent {\rm(b)}  If also $s_2\ne0$, then
  \begin{align}&\mhat{4}{a} {s_1, s_2, 
    s_3} =\frac{\biggl[\ds{ C_a(s_1,s_2) \mhat {4}{a} {s_1 + 1, s_2, 
    s_3}\atop-(1+s_1+s_2-s_3)\mhat {4}{a}{s_1 + 1, s_2+1, s_3 }}\biggr]  }{2s_2\prod_{k=1}^4(s_1+a_k)},  \label{prop1a} \end{align}where 
    \begin{align} &C_a(p,q)\nn\\&=q (p^2 + (p+q)^2)- 
\  q (a_1^2+a_2^2+a_3^2+a_1a_2 +a_1a_3+a_2a_3)\nn\\&  +(a_1 + a_2) (a_1 + a_3) (a_2 + a_3)\nn\\&=q (p^2 + (p+q)^2)- 
\frac{1}{2}  q\,(a_1^2 + a_2^2+ a_3^2+a_4^2)   -\frac{1}{3}(a_1^3 + a_2^3 +a_3^3+ a_4^3 ) \nn\\&=2  q(p+ a_1) (p+q - a_1) +  (q+a_2 + 
    a_3 )(q + a_2 + a_4)(q + a_3 + a_4) 
 .\label{Cadef}\end{align}  
 
\end{prop}

\begin{rem} {\rm The equivalence of the above three expressions for $C_a(p,q)$ is readily checked using straightforward algebra.  The first two of these expressions highlight the invariance of $C_a(p,q)$ under permutations of the $a_k$'s; the third expression will be of use in simplifying some calculations in Sections 3 and 4 below.}\end{rem}

\begin{pf} We first prove \eqref{prop1}. We multiply equation \eqref{gl4rec} through by $(s_{1}-s_2-a_1)(s_2-s_{3}-a_1)$.  From the resulting new equation, we subtract the same equation but with $a_1$ and $a_2$ interchanged. Using the fact that $\mhat{4}{a}{s}$ is invariant under permutations of the $a_j$'s, we thereby eliminate $\mhat{4}{a}{s_1,s_2+1,s_3}$. We find, after dividing everything through by $a_1-a_2$, that 
 \begin{align}  &(-s_1+s_3+a_1+a_2)T_{4,a}(s_1,s_2,s_3)\nn\\&+\frac{(s_1+s_2-s_3)T_{4,a}(s_1+1,s_2,s_3)}{( s_1+a_1)( s_{1}+a_2)} +\frac{(s_1-s_2-s_3)T_{4,a}(s_1,s_2,s_3+1)}{( s_3-a_1)( s_{3}-a_2)}\nn\\&-\frac{(s_1-s_3+a_1+a_2)T_{4,a}(s_1+1,s_2,s_3+1)}{( s_1+a_1)( s_{1}+a_2)(s_3-a_1)( s_{3}-a_2)}  =0.\label{gl4rec2}\end{align}  
 We now multiply \eqref{gl4rec2} through by $s_3-a_2$ and then subtract, from the resulting new equation, the same equation but with $a_2$ and $a_3$ interchanged.  This eliminates the term  $\mhat{4}{a}{s_1,s_2,s_3+1}$ and gives, upon division through by $(s_1+a_4) (a_2-a_3)$, the equation \begin{align}  T_{4,a}(s_1,s_2,s_3)&-\frac{(s_1+s_2-s_3)(s_1+s_3)T_{4,a}(s_1+1,s_2,s_3)}{( s_1+a_1)( s_{1}+a_2)( s_{1}+a_3)( s_{1}+a_4)} \nn\\&-\frac{ T_{4,a}(s_1+1,s_2,s_3+1)}{( s_1+a_1)( s_{1}+a_2)(s_1+a_3)( s_{1}+a_4)} =0, \end{align}  
which is the desired result.

We next prove (b). We eliminate $T_{4,a}(s_1+1,s_2,s_3)$ and $T_{4,a}(s_1+1,s_2,s_3+1)$ from \eqref{gl4rec}, much in the same way as we eliminated $T_{4,a}(s_1,s_2+1,s_3)$ and $T_{4,a}(s_1,s_2,s_3+1)$  above.  We obtain the following relation:
 \begin{align}  C_a(-s_3,s_2)T_{4,a}(s_1,s_2,s_3)&- (s_1+s_2-s_3) T_{4,a}(s_1,s_2+1,s_3) \nn\\& -2s_2 T_{4,a}(s_1,s_2,s_3+1)  =0.\label{gl4rec4}\end{align}   Into \eqref{gl4rec4}, we substitute $s_1\to s_1+1$, yielding a relation among $T_{4,a}(s_1+1,s_2,s_3)$, $T_{4,a}(s_1+1,s_2+1,s_3)$, and $T_{4,a}(s_1+1,s_2,s_3+1)$.  We may combine this latter relation with \eqref{prop1} to eliminate $T_{4,a}(s_1+1,s_2,s_3+1)$; the result is \begin{align}  &2 s_2 (s_1+a_1)(s_1+a_2)(s_1+a_3)(s_1+a_4)T_{4,a}(s_1,s_2,s_3)\nn\\&- C_a(s_1,s_2) T_{4,a}(s_1+1,s_2,s_3)  +(1+s_1+s_2-s_3) T_{4,a}(s_1+1,s_2+1,s_3)  =0,\label{gl4rec5}\end{align}   as desired.\hfill  $\Box$\end{pf} 
We now derive a relation entailing strictly positive shifts in $s_2$.
\begin{prop}\label{secondrecur} If $s_2\ne -a_j-a_k$  for $1\le j< k\le 4$, then
  \begin{align} &\mhat{4}{a} {s_1, s_2, 
    s_3}  =\frac{ \ds\biggl[{2 s_2 ( 1+s_1 + s_2 - s_3)\mhat {4}{a}{s_1 + 1, s_2+1, s_3 } \qquad\atop\qquad + ( s_2 + s_3-s_1 )C_a(s_1,s_2) \mhat {4}{a} {s_1 , s_2+1, 
    s_3} }\biggr]}{ \prod_{1\le j<k\le 4}( s_2+a_j+a_k) }.  \label{prop2} \end{align}   \end{prop}

\begin{pf}  We proceed as in the proofs of parts (a) and (b) of Proposition \ref{firstrecur}, but this time we eliminate, from \eqref{gl4rec}, the transforms $\mhat {4}{a} {s_1  , s_2, 
    s_3+1}$ and $\mhat {4}{a} {s_1+1 , s_2, 
    s_3+1}$.  The result is the relation    
    \begin{align}  C_a(s_1,-s_2) T_{4,a}(s_1 , s_2, s_3)& + 2 s_2 T_{4,a}(s_1+1, s_2, s_3 )\nn\\&+( s_2 +s_3-s_1) T_{4,a}( s_1, s_2 + 1, s_3) =0.\label{prop2first}\end{align}   We form a linear combination of  \eqref{prop2first} and \eqref{prop1a} to eliminate  $T_{4,a}(s_1+1 , s_2, s_3)$; the result is \eqref{prop2}. \hfill  $\Box$\end{pf} 
  We now deduce some relations entailing strictly positive shifts in $s_3$.  
\begin{prop}   \label{thirdrecur}
\noindent{\rm(a)}  If $s_3\ne a_k$ for $1\le k\le 4$, then 
\begin{equation}\mhat{4}{a} {s_1, s_2, 
    s_3} =\frac{ B(s_3,s_2,s_1) \mhat {4}{a} {s_1  , s_2, 
    s_3+1}+\mhat {4}{a}{s_1 + 1, s_2, s_3 + 1}  }{\prod_{k=1}^4(s_3-a_k)}.  \label{prop3a}    \end{equation} 
\noindent \rm(b) If $s_3\ne a_k$ for $1\le k\le 4$ and $s_2\ne0$, then
  \begin{align}&\mhat{4}{a} {s_1, s_2, 
    s_3} =\frac{\ds \biggl[{C_{-a}(s_3,s_2) \mhat {4}{a} {s_1 , s_2, 
    s_3+1}\atop-(1+ s_2+s_3-s_1)\mhat {4}{a}{s_1 , s_2+1, s_3+1 }}\biggr]  }{2s_2\prod_{k=1}^4(s_3-a_k)}.  \label{prop3b} \end{align}  
    
\noindent{\rm(c)}  If $s_2\ne -a_j-a_k$  for $1\le j< k\le 4$, then
  \begin{align} &\mhat{4}{a} {s_1, s_2, 
    s_3}  =\frac{\biggl[\ds{ 2 s_2 ( 1+s_2 + s_3 - s_1)\mhat {4}{a}{s_1 , s_2+1, s_3+1 } \qquad\atop\qquad+ ( s_1 + s_2-s_3 )C_{-a}(s_3,s_2) \mhat {4}{a} {s_1 , s_2+1, 
    s_3}}\biggr] }{ \prod_{1\le j<k\le 4}( s_2+a_j+a_k) }.  \label{prop3c} \end{align}   
    
 \end{prop}
        
\begin{pf}     It is well-known -- and may, in fact, be deduced from \eqref{Tnas}, and from induction -- that the Mellin transform $T_{n,a}(s)$ is invariant under the transformation\begin{align}&(s_1,s_2,\ldots ,s_{n-2},s_{n-1},a_1,a_2,\ldots, a_n) \nn\\&\rightarrow (s_{n-1},s_{n-2},\ldots , s_2,s_{1},-a_1,-a_2,\ldots, -a_n).\label{transf}\end{align} Parts (a), (b), and (c) of Proposition \ref{thirdrecur} then follow from Proposition \ref{firstrecur}(a), Proposition \ref{firstrecur}(b), and Proposition \ref{secondrecur} respectively. \hfill  $\Box$\end{pf} 
    For some applications, it is useful to have a recurrence relation expressing $T_{4,a}(s)$ in terms of ``strictly positive shifts in all $s_k$'s'' -- that is, in terms of transforms $T_{4,a}(s+\mu)$ where $\mu\in(\Z^+)^3$.  We conclude this section with such a relation.
\begin{prop} \label{fourthrecur} Let $B(p,q,r)$ be as in \eqref{Bdef}, and $C_a(p,q)$   as in \eqref{Cadef}.   If $s_1\ne -a_k$ and $s_3\ne a_k$ for $1\le k\le 4$, and $s_2\ne -a_j-a_k$  for $1\le j< k\le4$, then
    \begin{align} &\mhat{4}{a} {s_1, s_2, 
    s_3}=\biggl[{\prod_{k=1}^4 (s_1+a_k) (s_3-a_k)\prod_{ j=1}^{k-1}( s_2+a_j+a_k)}\biggr]^{-1}\nn
    \\&\times\biggl\{
 \biggl[2 s_2\prod_{k=1}^4 (s_1+a_k)+  B(s_3,s_2 ,s_1) C_a(s_1,s_2)\biggr]( 1+s_1 + s_2 - s_3)\nn\\&\times  \mhat {4}{a}{s_1 + 2, s_2+1, s_3+1 } + \biggl(
  \biggl[2 s_2 \prod_{k=1}^4 (s_1+a_k) +B( s_3 ,s_2,s_1 )C_a(s_1,s_2)  \biggr]\nn\\&\times B( s_1+1,s_2,s_3)  
+C_a(s_1,s_2)\prod_{k=1}^4  (s_3-a_k) \biggr)  ( s_2 + s_3-s_1 )\nn\\&\times
 \mhat {4}{a} {s_1+1 , s_2+1, 
    s_3+1}\biggr\} . \end{align}     \end{prop}
\begin{pf}  We substitute $s_2\to s_2+1$ into  \eqref{prop1}, to get   \begin{align}&\mhat{4}{a} {s_1, s_2+1, 
    s_3}\nn\\& =\frac{ B(s_1,s_2+1,s_3) \mhat {4}{a} {s_1 + 1, s_2+1, 
    s_3}+\mhat {4}{a}{s_1 + 1, s_2+1, s_3 + 1}  }{\prod_{k=1}^4(s_1+a_k)}.  \label{propb} \end{align}We then apply this result to \eqref{prop2} to eliminate $\mhat{4}{a} {s_1, s_2+1, 
    s_3}$; we get  \begin{align} &\mhat{4}{a} {s_1, s_2, 
    s_3}\nn  \\&  =\biggl[{\prod_{k=1}^4 (s_1+a_k) \prod_{ j=1}^{k-1}( s_2+a_j+a_k)}\biggr]^{-1}
 \biggl\{
\biggl(2 s_2 ( 1+s_1 + s_2 - s_3) \prod_{k=1}^4 (s_1+a_k)\nn\\&+ ( s_2 + s_3-s_1 )C_a(s_1,s_2)B(s_1 , s_2+1, s_3)   \biggr)\mhat {4}{a}{s_1 + 1, s_2+1, s_3 } \nn\\&+ ( s_2 + s_3-s_1 )C_a(s_1,s_2) \mhat {4}{a} {s_1+1 , s_2+1, 
    s_3+1}\biggr\} .  \label{prop1c} \end{align}  
    Next, we recall that $\mhat{4}{a}{s}$ is invariant under the transformation \eqref{transf}, so that \eqref{propb} yields\begin{align}&\mhat{4}{a} {s_1, s_2+1, 
    s_3}\nn\\& =\frac{ B(s_3,s_2+1,s_1) \mhat {4}{a} {s_1 , s_2+1, 
    s_3+1}+\mhat {4}{a}{s_1 + 1, s_2+1, s_3 + 1}  }{\prod_{k=1}^4(s_3-a_k)}.  \label{propish} \end{align}Substituting $s_1\to s_1+1$ into this then gives  
    \begin{align} &\mhat{4}{a} {s_1+1, s_2+1, 
    s_3} \nn\\&\hskip-1pt=\frac{ B(s_3,s_2+1,s_1+1) \mhat {4}{a} {s_1+1 , s_2+1, 
    s_3+1}+\mhat {4}{a}{s_1 + 2, s_2+1, s_3 + 1}  }{\prod_{k=1}^4(s_3-a_k)}.  \label{prop4d} \end{align}   We now put \eqref{prop4d} into \eqref{prop1c} to eliminate $\mhat{4}{a}{s_1+1,s_2+1,s_3}$; the result is
    \begin{align} &\mhat{4}{a} {s_1, s_2, 
    s_3}  =\biggl[{\prod_{k=1}^4 (s_1+a_k) (s_3-a_k)\prod_{ j=1}^{k-1}( s_2+a_j+a_k)}\biggr]^{-1}
\nn\\&\times  \biggl\{
 \biggl(2 s_2( 1+s_1 + s_2 - s_3) \prod_{k=1}^4 (s_1+a_k)+  B(s_1,s_2+1,s_3)( s_2 + s_3-s_1 )\nn\\&\times C_a(s_1,s_2)\biggr)\mhat {4}{a}{s_1 + 2, s_2+1, s_3+1 }+ \biggl(( s_2 + s_3-s_1 )C_a(s_1,s_2)\nn\\&\times \prod_{k=1}^4  (s_3-a_k) +
B( s_3, s_2+1 ,s_1+1) 
 \biggl[2 s_2 (s_1+s_2-s_3+1)\prod_{k=1}^4 (s_1+a_k) \nn \\& +B( s_1 ,s_2+1,s_3 )(s_2+s_3-s_1)C_a(s_1,s_2)  \biggr]\biggr)
 \mhat {4}{a} {s_1+1 , s_2+1, 
    s_3+1}\biggr\} .  \end{align}  Some rearrangement of these terms, using, for example, the fact that\begin{align*} &B(s_1,s_2+1,s_3)( s_2 + s_3-s_1 )= (1+s_1+s_2-s_3)(s_1+s_3)( s_2 + s_3-s_1 )\\&=(1+s_1+s_2-s_3)B(s_3,s_2,s_1),\end{align*}   then yields the stated result. \hfill  $\Box$\end{pf} 

\section{Poles and residues of \boldmath ${T_{n,a}(s)}$}

As described in Section 1 above, recent work \cite {GW}  on orthogonality of Maass forms for $GL(4)$ involves a test function $p_{T,R}$ that is defined as an integral of a known function  $p_{T,R}^\sharp$ (the Lebedev-Whittaker transform of $p_{T,R}$) against a Whittaker function.  

To obtain estimates of the desired strength for $p_{T,R}$, it is necessary to move certain lines of integration, in this integral defining $p_{T,R}$, a finite distance to the left, and to analyze the poles and residues of $T_{4,a}(s)$ that are thereby encountered.  We develop the requisite results concerning these poles and residues here.

From the analyticity of $T_{4,a}(s)$ for $s_k$'s of sufficiently large real parts, and from Proposition \ref{fourthrecur}, we may readily deduce the following.

\begin{prop}\label{pole-loc}The Mellin transform $T_{4,a}(s)$ extends to a meromorphic function of the variable $(s_1,s_2,s_3)\in\C^3$,   with poles at
\begin{alignat*}{2}s_1&=\phantom{-a_\ell\ }-a_m-\delta_1\qquad&(m\in\{1,2,3,4\},\ \delta_1\in\Z_{\ge0});
\\ s_2&=-a_m-a_n-\delta_2\qquad&(\{m,n\}\subset\{1,2,3,4\},\  \delta_2\in\Z_{\ge0});\\
s_3&=\phantom{-a_\ell-\ \ }a_m-\delta_3\qquad&(m\in\{1,2,3,4\},\ \delta_3\in\Z_{\ge0}),\end{alignat*}and no other poles or polar divisors in $\C^3$.\end{prop}

The analyses of \cite{GW} require explicit information concerning the residues of $T_{4,a}(s)$ at the above poles.  Using our above recurrence relations, we will develop some formulas, below,  for these residues.  These formulas will entail a certain polynomial  $p_\delta$, the salient properties of which we now describe.

\begin{lem} \label{pdelta-lem} For $a\in\C$, define\begin{equation}(a)_\kappa=(a+\kappa-1)(a+\kappa-2)\cdots
(a+1)a\quad\hbox{$(\kappa\in\Z^+)$;\qquad$(a)_0=1$}.\label{poch1} \end{equation}Then, for $\delta\in\Z_{\ge0}$ and $b,c,d,e,f,g\in\Z$, define the  polynomial \begin{equation}p_\delta(b,c,d;e,f,g)=p_\delta\biggl[{b,c,d;\atop e,f,g}\biggr],\label{pdelta}\end{equation}  of degree $\le3\delta$,  by
\begin{align}p_\delta(b,c,d;e,f,g) &=\sum_{\kappa=0}^\delta(-1)^\kappa\frac{ (b)_\kappa (c)_\kappa(d)_\kappa (e+\kappa)_{\delta-\kappa}(f+\kappa)_{\delta-\kappa} (g+k )_{\delta-\kappa}}{\kappa!(\delta-\kappa)!}.\label{pdef}\end{align} 
 \noindent{\rm (a)} We have\begin{align}&(e+\delta)(f+\delta)(g -1 )p_\delta\biggl[{b,c,d;\atop e,f,g}\biggr]-b\,c\,d\, p_\delta\biggl[{b+1,c+1,d+1;\atop e+1,f+1,g+1}\biggr]\nn\\&=(\delta+1) \,p_{\delta+1}\biggl[{b,c,d;\atop  e,f,g -1 }\biggr].\label{pdelta-recur1}\end{align}
    \noindent{\rm (b)} Suppose $e+f+g+\delta-b-c-d=1$.  Then  \begin{align}    & b\,  (f - c) (f-d)(c-e-(1+ \delta )) \,p_\delta\biggl[{ b +1, c, 
   d  ;\atop e, f+1, 
   g  }\biggr]  \nr&+ (g -b-1)\bigl((e-b+\delta)(e-c)  (f-1) -b\,(f-c)(1+\delta)\bigr)   p_\delta\biggl[{ b , c, d;\atop e, f,
    g}\biggr]\nr& =(\delta + 1) (  e-c ) p_{\delta + 1}\biggl[{b  , c-1, d;\atop e,
    f-1, g-1}\biggr] .\label{pdelta-recur2}
    \end{align}
  \noindent{\rm (c)} If any of the variables $b,c$, or $d$ equals a nonpositive integer $-\gamma$, where $0\le \gamma \le \delta$, then $p_\delta(b,c,d;e,f,g)$ is divisible by \begin{align*}(e+\gamma)_{\delta-\gamma}(f+\gamma)_{\delta-\gamma} (g+\gamma)_{\delta-\gamma} & =\prod_{\kappa=0}^{\delta-\gamma-1}(e+\gamma+\kappa)(f+\gamma+\kappa) (g+\gamma+\kappa)  \\
  & = \frac{\gg{e+\delta}\gg{f+\delta} \gg{g+\delta}}{\gg{e+\gamma}\gg{f+\gamma}\gg{g+\gamma} }   .\end{align*}
\end{lem}
\begin{pf} We begin with part (a).  By definition of $p_\delta$, we see that   \begin{align}& (e+\delta)(f+\delta)(g-1 )p_\delta\biggl[{b,c,d;\atop e,f,g}\biggr]-b\,c\,d\, p_\delta\biggl[{b+1,c+1,d+1;\atop e+1,f+1,g+1}\biggr]
\nn\\&=
  (e+\delta)(f+\delta)(g-1 ) 
 \sum_{\kappa=0}^\delta(-1)^\kappa  \frac{(b)_\kappa(c)_\kappa (d)_\kappa
 (e+\kappa)_{\delta-\kappa} (f+\kappa)_{\delta-\kappa} (g+\kappa )_{\delta-\kappa}}{\kappa!(\delta-\kappa)!}
 \nr&+ b\,c\,d\, \sum_{\kappa=0}^\delta (-1)^\kappa\frac{ (b+1)_\kappa (c+1)_\kappa (d+1)_\kappa
}{\kappa!(\delta-\kappa)!} \nr&\times  (e+1+\kappa)_{\delta-\kappa}(f+1+\kappa)_{\delta-\kappa}(g+1+\kappa )_{\delta-\kappa}
 \nn\\&=     \sum_{\kappa=0}^\delta (-1)^\kappa\frac{ (b)_\kappa(c)_\kappa (d)_\kappa
 (e+\kappa)_{\delta+1-\kappa}  (f+\kappa)_{\delta+1-\kappa}(g-1 )(g +\kappa)_{\delta -\kappa}}{\kappa!(\delta-\kappa)!}
\nn\\&-\hskip-3pt\sum_{\kappa=0}^{\delta}(-1)^\kappa \frac{ (b)_{\kappa+1} (c)_{\kappa+1} (d)_{\kappa+1}
 (e+1+\kappa)_{\delta-\kappa}(f+1+\kappa)_{\delta-\kappa}(g +1 +\kappa )_{\delta-\kappa}}{\kappa!(\delta-\kappa)!} ,\label{pdelta-a}\end{align}
 the last step because 
\begin{align*} (e+\delta)(f+\delta) (e+\kappa)_{\delta-\kappa} (f+\kappa)_{\delta-\kappa}  &=(e+\kappa)_{\delta+1-\kappa} (f+\kappa)_{\delta+1-\kappa}\end{align*}and\begin{align*}b\,c\,d\, (b+1)_\kappa (c+1)_\kappa (d+1)_\kappa& = (b)_{\kappa+1} (c)_{\kappa+1} (d)_{\kappa+1} .\end{align*}We substitute $\kappa\to \kappa-1$ into the second sum on the right-hand side of \eqref{pdelta-a}, to get
 \begin{align}&(e+\delta)(f+\delta)(g-1 )p_\delta\biggl[{b,c,d;\atop e,f,g}\biggr]-b\,c\,d\, p_\delta\biggl[{b+1,c+1,d+1;\atop e+1,f+1,g+1}\biggr]\nn\\& =\hskip-3pt \sum_{\kappa=0}^\delta (-1)^\kappa\frac{(b)_\kappa(c)_\kappa (d)_\kappa
 (e+\kappa)_{\delta+1-\kappa}  (f+\kappa)_{\delta+1-\kappa} (g-1 )(g+\kappa)_{\delta-\kappa}}{\kappa!(\delta-\kappa)!}
\nn\\&- \sum_{\kappa=1}^{\delta+1} (-1)^{\kappa-1}\frac{(b)_{\kappa} (c)_{\kappa} (d)_{\kappa}
 (e+\kappa)_{\delta+1-\kappa}(f+\kappa)_{\delta+1-\kappa}(g +\kappa )_{\delta+1-\kappa}}{(\kappa-1)!(\delta+1-\kappa)!} 
 \nn\\& = \frac{ 
 (e)_{\delta+1}  (f)_{\delta+1} (g-1 ) (g )_{\delta}}{ \delta!} +(-1)^{\delta+1} \frac{(b)_{\delta+1} (c)_{\delta+1} (d)_{\delta+1}
  }{\delta! } \nn\\&
 +  \sum_{\kappa=1}^\delta(-1)^\kappa\biggl(\frac{(g-1 )(g+\kappa)_{\delta-\kappa}}{\kappa!(\delta-\kappa)!}+\frac{ (g +\kappa )_{\delta+1-\kappa}}{(\kappa-1)!(\delta+1-\kappa)!}\biggr)\nn\\&\
 \times{(b)_\kappa(c)_\kappa (d)_\kappa
 (e+\kappa)_{\delta+1-\kappa}  (f+\kappa)_{\delta+1-\kappa}} 
 ,\label{pdelta-b}\end{align}  the last step by separating the $\kappa=0$ term off from the first sum,  separating the $\kappa=\delta+1$ term off from the second sum, and then combining the remaining two sums on $\kappa$ (from $\kappa=1$ to $\kappa=\delta$).  But $$(g-1 ) (g )_{\delta} =(g-1 )_{\delta+1},$$and\begin{align*}&\frac{(g-1 )(g+\kappa)_{\delta-\kappa}}{\kappa!(\delta-\kappa)!}+\frac{ (g +\kappa )_{\delta+1-\kappa}}{(\kappa-1)!(\delta+1-\kappa)!}
  \\&=\frac{(g+\kappa)_{\delta-\kappa}}{\kappa!(\delta+1-\kappa)!} \bigl({(\delta+1-\kappa)(g-1 ) } +{\kappa} (g+\delta)\bigr)
   \\&=\frac{(g+\kappa)_{\delta-\kappa}}{\kappa!(\delta+1-\kappa)!} {(\delta+1) (g-1 +\kappa)}
    =(\delta+1) \frac{(g-1 +\kappa)_{\delta+1-\kappa}}{\kappa!(\delta+1-\kappa)!}   .\end{align*}
  So \eqref{pdelta-b} gives 
  \begin{align*}&(e+\delta)(f+\delta)(g-1 )p_\delta\biggl[{b,c,d;\atop e,f,g}\biggr]-b\,c\,d\, p_\delta\biggl[{b+1,c+1,d+1;\atop e+1,f+1,g+1}\biggr]\nn\\&  
 = (\delta+1) \biggl[   \frac{ 
 (e)_{\delta+1}  (f)_{\delta+1} (g-1  )_{\delta+1}}{ (\delta+1)!}+(-1)^{\delta+1} \frac{(b)_{\delta+1} (c)_{\delta+1} (d)_{\delta+1}
  }{(\delta+1)! }
\nn\\& 
 +  \sum_{\kappa=1}^\delta (-1)^\kappa\frac{ (b)_\kappa(c)_\kappa (d)_\kappa
 (e+\kappa)_{\delta+1-\kappa}  (f+\kappa)_{\delta+1-\kappa} (g-1 +\kappa)_{\delta+1-\kappa}}{\kappa!(\delta+1-\kappa)!}\biggr] 
  \nn\\& =  (\delta+1) 
  \sum_{\kappa=0}^{\delta +1} (-1)^\kappa \frac{ (b)_\kappa(c)_\kappa (d)_\kappa
 (e+\kappa)_{\delta+1-\kappa}  (f+\kappa)_{\delta+1-\kappa} (g-1 +\kappa)_{\delta+1-\kappa}}{\kappa!(\delta+1-\kappa)!}
   \nn\\&  
 =(\delta+1)\, p_{\delta+1}\biggl[{b,c,d;\atop e,f,g-1 }\biggr], \end{align*}  and part (a) of our lemma is proved.
 
 The proof of part (b) is similar (though somewhat messier). We omit the details.

 To prove part (c), we assume that  $b=-\gamma$ where $\gamma$ is an integer and $0\le \gamma\le \delta$.  Then, by \eqref{poch1}, we see that $(b)_\kappa=0$ for $\kappa>\gamma$, so \eqref{pdef} gives
\begin{align}&p_\delta(b,c,d;e,f,g) =\sum_{\kappa=0}^\gamma(-1)^\kappa\frac{ (b)_\kappa (c)_\kappa(d)_\kappa (e+\kappa)_{\delta-\kappa}(f+\kappa)_{\delta-\kappa} (g+\kappa )_{\delta-\kappa}}{\kappa!(\delta-\kappa)!}\nn
\\&=(e+\gamma)_{\delta-\gamma}(f+\gamma)_{\delta-\gamma}(g+\gamma)_{\delta-\gamma}
\nr&\times\sum_{\kappa=0}^\gamma(-1)^\kappa\frac{ (b)_\kappa (c)_\kappa(d)_\kappa (e+\kappa)_{\gamma-\kappa}(f+\kappa)_{\gamma-\kappa} (g+\kappa )_{\gamma-\kappa}}{\kappa!(\delta-\kappa)!},\label{pdelta-short}\end{align} 
the last step because, by \eqref{poch1},  $$(a+\kappa)_{\delta-\kappa}=(a+\gamma)_{\delta-\gamma}(a+\kappa)_{\gamma-\kappa}$$for $\kappa\le\gamma\le \delta$. Certainly  \eqref{pdelta-short} implies that $(e+\gamma)_{\delta-\gamma}(f+\gamma)_{\delta-\gamma}(g+\gamma)_{\delta-\gamma}$ divides $p_\delta(b,c,d;e,f,g)$, so part (b) of our lemma is proved.  \hfill  $\Box$\end{pf}
  Using the above lemma, we may now deduce some explicit results concerning the residues of $T_{4,a}(s)$ at the poles described in Proposition \ref{pole-loc}.
\begin{prop}   \label{singleres}
\noindent{\rm(a)}  For  $m \in \{1,2,3,4\} ,$  $\{n,p,q\}= \{1,2,3,4\}\backslash \{m\}$, and $\delta_1\in\Z_{\ge0}$, we have\begin{align}&\res{s_1=-a_m-\delta_1}{\mhat{4}{a}{s}}\nn\\&=\frac{1}{ \gg{s_2+s_3+a_m+\delta_1}}\biggl[\prod_{k\in\{n,p,q\} }\gg{a_k-a_m-\delta_1}\gg{s_3-a_k}
\gg{s_2+a_k+a_m}\biggr]\nn
\\&\times p_{\delta_1} (s_3-a_n,s_3-a_p,s_3-a_q;s_2+s_3+a_m, 1+a_m-s_2+s_3,s_3-a_m-\delta_1).\label{s1-res}\end{align}
\noindent{\rm(b)} For  $\{m,n\} \subset\{1,2,3,4\} ,$  $\{p,q\}= \{1,2,3,4\}\backslash \{m,n\}$, and $\delta_2\in\Z_{\ge0}$, we have\begin{align}&\res{s_2=-a_m-a_n-\delta_2}{\mhat{4}{a}{s}}\nn
    \\& = \biggl[\prod_{j\in\{m,n\} }\gg{s_1+a_j}\biggr]\biggl[\prod_{ k\in\{p,q\}}\gg{s_3-a_k} \prod_{j\in\{m,n\}  }\gg{a_k-a_j-\delta_2}\biggr]    \nn     \\& \times p_{\delta_2}({s_1+a_m, a_p-a_n-\delta_2,s_3-a_q;  1+a_m-a_q,s_1+a_p-\delta_2,s_3-a_n-\delta_2}).\label{s2-res}\end{align}
\noindent{\rm (c)} For  $m \in \{1,2,3,4\} ,$  $\{n,p,q\}= \{1,2,3,4\}\backslash \{m\}$, and $\delta_3\in\Z_{\ge0}$, we have \begin{align*}&\res{s_3=a_m-\delta_3}{\mhat{4}{a}{s}}\nn\\&=\frac{1}{ \gg{s_1+s_2-a_m+\delta_3}}\biggl[\prod_{k\in\{n,p,q\} }\gg{a_m-a_k-\delta_3}\gg{s_1+a_k}
\gg{s_2-a_k-a_m}\biggr]\nn
\\&\times p_{\delta_3} (s_1+a_n,s_1+a_p,s_1+a_q;s_1+s_2-a_m, 1-a_m-s_2+s_1,s_1+a_m-\delta_3). \end{align*}\end{prop}

\begin{pf}  We begin with   part (a) of our proposition.  We will prove the desired result by induction on $\delta_1$.

The case  $\delta_1=0$ is given given by \cite[Theorem 3.2]{ST1}, in the case $n=4$.  So let us now assume that \eqref{s1-res} is true for a nonnegative integer $\delta_1=\delta$.  From the recurrence relation \eqref{prop1}, we have
  \begin{align*} &\hbox{\rm Res}_{s_1=-a_m-(\delta +1)}\mhat{4}{a} {s_1, s_2, 
    s_3}\nn\\ &=\hbox{\rm Res}_{s_1=-a_m-(\delta +1)}\biggl(\frac{ B(s_1,s_2,s_3) \mhat {4}{a} {s_1 + 1, s_2, 
    s_3}+\mhat {4}{a}{s_1 + 1, s_2, s_3 + 1}  }{\prod_{k=1}^4(s_1+a_k)}\biggr) 
    \nn\\ &=\hbox{\rm Res}_{s_1=-a_m-\delta }\biggl(\frac{ B(s_1-1,s_2,s_3) \mhat {4}{a} {s_1 , s_2, 
    s_3}+\mhat {4}{a}{s_1 , s_2, s_3 + 1}  }{\prod_{k=1}^4(s_1-1+a_k)}\biggr) ,
   \nn \end{align*}so that, by the induction hypothesis,\begin{align*}&\hbox{\rm Res}_{s_1=-a_m-(\delta +1)}\mhat{4}{a} {s_1, s_2, 
    s_3}=\frac{1}{-(\delta+1)}
   \nr&\times\biggl(  \frac{ B(-a_m-\delta-1,s_2,s_3)}{ \gg{s_2+s_3+a_m+\delta }}
  \biggl[\prod_{k\in\{n,p,q\} }\hskip-9pt\frac{\gg{a_k-a_m-\delta }\gg{s_3-a_k}
\gg{s_2+a_k+a_m}}{a_k-a_m-\delta-1}\biggr] \nn\\&\times p_{\delta} \biggl[{s_3-a_n,s_3-a_p,s_3-a_q;\atop s_2+s_3+a_m, 1+a_m-s_2+s_3,s_3-a_m-\delta}\biggr] +
\frac{1}{ \gg{s_2+s_3+1+a_m+\delta}}\nn\\&\times\biggl[\prod_{k\in\{n,p,q\} }\hskip-9pt\frac{\gg{a_k-a_m-\delta }\gg{s_3+1-a_k}
\gg{s_2+a_k+a_m}}{a_k-a_m-\delta-1}\biggr]\nn
\\&\times p_{\delta} \biggl[{s_3+1-a_n,s_3+1-a_p,s_3+1-a_q;\atop s_2+s_3+1+a_m, 2+a_m-s_2+s_3,s_3+1-a_m-\delta}\biggr]\biggr). 
  \end{align*}Using the translation equation $\gg{s+1}=s\gg{s}$, and the definition \eqref{Bdef} of $B(s_1,s_2,s_3)$, we may rewrite the above as follows:
    \begin{align} &\hbox{\rm Res}_{s_1=-a_m-(\delta +1)}\mhat{4}{a} {s_1, s_2, 
    s_3} = \frac{1}{-(\delta+1) \gg{s_2+s_3+a_m+\delta+1}}
 \nn\\ &\times    \biggl[\prod_{k\in\{n,p,q\} }{\gg{a_k-a_m-\delta -1}\gg{s_3-a_k}
\gg{s_2+a_k+a_m}}\biggr] \nn\\&\times\biggl(   { (-a_m-\delta-1+s_2-s_3)(-a_m-\delta-1+s_3)(s_2+s_3+a_m+\delta )} \nn\\&\times  p_{\delta} \biggl[{s_3-a_n,s_3-a_p,s_3-a_q;\atop s_2+s_3+a_m, 1+a_m-s_2+s_3,s_3-a_m-\delta}\biggr]\nn\\&+
 (s_3-a_n)(s_3-a_p)(s_3-a_q)\nn
\\&\times p_{\delta} \biggl[{s_3+1-a_n,s_3+1-a_p,s_3+1-a_q;\atop s_2+s_3+1+a_m, 2+a_m-s_2+s_3,s_3+1-a_m-\delta }\biggr]\biggr). \label{T4-recur-a}
  \end{align}But by Lemma \ref{pdelta-lem}(a), the quantity in large parentheses, in \eqref{T4-recur-a}, equals$$-(\delta+1)p_{\delta+1} \biggl[{s_3-a_n,s_3-a_p,s_3-a_q;\atop s_2+s_3+a_m, 1+a_m-s_2+s_3,s_3-a_m-\delta-1}\biggr].$$So \eqref{T4-recur-a} yields
    \begin{align*} &\hbox{\rm Res}_{s_1=-a_m-(\delta +1)}\mhat{4}{a} {s_1, s_2, 
    s_3} = \frac{1}{ \gg{s_2+s_3+a_m+\delta+1}}
 \nn\\ &\times    \biggl[\prod_{k\in\{n,p,q\} }{\gg{a_k-a_m-\delta -1}\gg{s_3-a_k}
\gg{s_2+a_k+a_m}}\biggr] \nn\\&\times p_{\delta+1} \biggl[{s_3-a_n,s_3-a_p,s_3-a_q;\atop s_2+s_3+a_m, 1+a_m-s_2+s_3,s_3-a_m-(\delta+1)}\biggr] ,    \end{align*}which tells us that \eqref{s1-res} is true for $\delta_1=\delta+1$ as well.  So by induction, part (a) of our lemma is proved.

We now prove part (b) of our proposition, by induction on $\delta_2$. We first note that the anchor step -- the case $\delta_2=0$ -- is given by \cite[Theorem 3.2]{ST1}, in the case $n=4$.  Now assume that \eqref{s2-res} is true for a nonnegative integer $\delta_2=\delta$.  From the recurrence relation \eqref{prop2}, we have
  \begin{align*} &{\rm Res}_{s_2=-a_m-a_n-(\delta+1)}\mhat{4}{a} {s_1, s_2, 
    s_3}  \nn\\&={\rm Res}_{s_2=-a_m-a_n-(\delta+1)}  \frac{ \ds\biggl[{2 s_2 ( 1+s_1 + s_2 - s_3)\mhat {4}{a}{s_1 + 1, s_2+1, s_3 } \qquad\atop\qquad + ( s_2 + s_3-s_1 )C_a(s_1,s_2) \mhat {4}{a} {s_1 , s_2+1, 
    s_3} }\biggr]}{ \prod_{1\le j<k\le 4}( s_2+a_j+a_k) }
\nn\\&={\rm Res}_{s_2=-a_m-a_n-\delta}  \frac{ \ds\biggl[{2( s_2-1) (  s_1 + s_2 - s_3)\mhat {4}{a}{s_1 + 1, s_2,s_3 }  \atop + ( s_2 + s_3-s_1 -1)C_a(s_1,s_2-1) \mhat {4}{a} {s_1 , s_2 , 
    s_3} }\biggr]}{ \prod_{1\le j<k\le 4}( s_2-1+a_j+a_k) } 
    \nn\\&= \biggl[(\delta+1)(2a_m+2a_n+\delta+1)\prod_{j\in\{m,n\}} \prod_{ k\in\{p,q\}}( a_k-a_j-\delta-1)\biggr]^{-1}
    \\&\times  \ds\bigl\{2( -a_m-a_n-\delta-1) (  s_1   - s_3-a_m-a_n-\delta)\\&\times{\rm Res}_{s_2=-a_m-a_n-\delta} \mhat {4}{a}{s_1 + 1, s_2,s_3 } \\&+ (   s_3-s_1-a_m-a_n-\delta -1)C_a(s_1,-a_m-a_n-\delta-1)\\&\times {\rm Res}_{s_2=-a_m-a_n-\delta} \mhat {4}{a} {s_1 , s_2 , 
    s_3} \bigr\},\end{align*}so that, by the induction hypothesis,
 \begin{align*} &{\rm Res}_{s_2=-a_m-a_n-(\delta+1)}\mhat{4}{a} {s_1, s_2, 
    s_3} \\&= \biggl[(\delta+1)(2a_m+2a_n+\delta+1)\prod_{j\in\{m,n\}} \prod_{ k\in\{p,q\}}( a_k-a_j-\delta-1)\biggr]^{-1}
    \\&\times   \biggl\{ 2 (-a_m-a_n-\delta-1)(  s_1   - s_3-a_m-a_n-\delta)
    \nn\\&\times  \biggl[\prod_{j\in\{m,n\} }\gg{s_1+1+a_j}\biggr]\biggl[\prod_{ k\in\{p,q\}}\gg{s_3-a_k} \prod_{j\in\{m,n\}  }\gg{a_k-a_j-\delta }\biggr]    \nn     \\& \times p_{\delta }\biggl[{s_1+1+a_m, a_p-a_n-\delta ,s_3-a_q;\atop 1+a_m-a_q,s_1+1+a_p-\delta ,s_3-a_n-\delta }\biggr] \nn\\&+ (  s_3-s_1-a_m-a_n-\delta -1)C_a(s_1,-a_m-a_n-\delta-1) 
    \nn\\&\times\biggl[\prod_{j\in\{m,n\} }\gg{s_1+a_j}\biggr]\biggl[\prod_{ k\in\{p,q\}}\gg{s_3-a_k} \prod_{j\in\{m,n\}  }\gg{a_k-a_j-\delta }\biggr]    \nn     \\& \times p_{\delta }\biggl[{s_1+a_m, a_p-a_n-\delta ,s_3-a_q;\atop 1+a_m-a_q,s_1+a_p-\delta ,s_3-a_n-\delta }\biggr]\biggr\} . \end{align*} 
    Using the relation $\gg{s+1}=s\,\gg{s}$, we rewrite the above as follows:
      \begin{align}
       &{\rm Res}_{s_2=-a_m-a_n-(\delta+1)}\mhat{4}{a} {s_1, s_2, 
    s_3} = \bigl[(\delta+1)(2a_m+2a_n+\delta+1) \bigr]^{-1}   \nn\\&\times \biggl[\prod_{j\in\{m,n\} }\gg{s_1+a_j}\biggr]\biggl[\prod_{ k\in\{p,q\}}\gg{s_3-a_k} \prod_{j\in\{m,n\}  }\gg{a_k-a_j-(\delta+1) }\biggr] \nn\\&\times     \biggl\{ 2 (-a_m-a_n-\delta-1)(  s_1   - s_3-a_m-a_n-\delta)(s_1+a_m)(s_1+a_n)
   \nn     \\& \times  p_{\delta }\biggl[{s_1+1+a_m, a_p-a_n-\delta ,s_3-a_q;\atop 1+a_m-a_q,s_1+1+a_p-\delta ,s_3-a_n-\delta }\biggr] \nn\\&+ (  s_3-s_1-a_m-a_n-\delta -1)C_a(s_1,-a_m-a_n-\delta-1) 
    \nn     \\& \times  p_{\delta }\biggl[{s_1+a_m, a_p-a_n-\delta ,s_3-a_q;\atop 1+a_m-a_q,s_1+a_p-\delta ,s_3-a_n-\delta }\biggr]\biggr\} . \label{s2-penult}\end{align}If we now define$$\displaylines{ b=s_1+a_m,\quad c= a_p-a_n-\delta , \quad d=s_3-a_q,\cr e=1+a_m-a_q, \quad f=s_1+a_p -\delta,\quad g= s_3-a_n-\delta}$$(so that $e+f+g+\delta-b-c-d=1$), then we check that the quantity in large braces, in \eqref{s2-penult}, is equal to 
    \begin{align*}    & b\,  (f - c) (f-d)(c-e-(1+ \delta )) \,p_\delta\biggl[{ b +1, c, 
   d  ;\atop e, f+1, 
   g  }\biggr]  \nr&+ (g -b-1)\bigl((e-b+\delta)(e-c)  (f-1) -b\,(f-c)(1+\delta)\bigr)   p_\delta\biggl[{ b , c, d;\atop e, f,
    g}\biggr].\end{align*}But by Lemma \ref{pdelta-lem}(b), this quantity equals
    
    \begin{align*}&(\delta + 1) (  e-c ) p_{\delta + 1}\biggl[{b  , c-1, d;\atop e,
    f-1, g-1}\biggr]= (\delta + 1) (2a_m+2a_n+\delta+1)\\&\times p_{\delta + 1}\biggl[{s_1+a_m, a_p-a_n-(\delta +1), s_3-a_q;\atop  1+a_m-a_q, s_1+a_p -(\delta+1), s_3-a_n-(\delta+1)}\biggr].\end{align*}So \eqref{s2-penult} reads
      \begin{align}
       &{\rm Res}_{s_2=-a_m-a_n-(\delta+1)}\mhat{4}{a} {s_1, s_2, 
    s_3}    \nn\\&=\biggl[\prod_{j\in\{m,n\} }\gg{s_1+a_j}\biggr]\biggl[\prod_{ k\in\{p,q\}}\gg{s_3-a_k} \prod_{j\in\{m,n\}  }\gg{a_k-a_j-(\delta+1) }\biggr] \nn\\&\times  p_{\delta + 1}\biggl[{s_1+a_m, a_p-a_n-(\delta +1), s_3-a_q;\atop  1+a_m-a_q, s_1+a_p -(\delta+1), s_3-a_n-(\delta+1)}\biggr]    ,\label{s2-ult}\end{align}which is to say that \eqref{s2-res} holds for $\delta_2=\delta+1$.  So, by induction, part (b) of our proposition is proved.
 
Part (c) of our proposition follows from part (a), and the fact that $\mhat{4}{a}{s}$ is invariant under the transformation \eqref{transf}.  \hfill  $\Box$\end{pf} 

Note that, because $T_{4,a}(s)$ is invariant under permutations of the $a_k$'s, so are the results of the above proposition.

We now wish to compute residues of the above residues in either of the remaining variables.  For example, we wish to compute  
 $$\hbox{\rm Res}_{s_1=-a_1-\delta_1} \bigl(\hbox{\rm Res}_{s_2=-a_1-a_4-\delta_2}{\mhat{4}{a}{s}}\bigr).$$To reduce the total number of computations that we'll need, we make several observations, which are clear from general principles and from the above proposition.
 
 \begin{itemize}
 \item ``Order doesn't matter;'' for example,
  \begin{align*}&\hbox{\rm Res}_{s_1=-a_1-\delta_1}\bigl(\hbox{\rm Res}_{s_3=a_2-\delta_3} {\mhat{4}{a}{s}}\bigr)=\hbox{\rm Res}_{s_3=a_2-\delta_3} \bigl(\hbox{\rm Res}_{s_1=-a_1-\delta_1}{\mhat{4}{a}{s}}\bigr).\end{align*}

  \item  The residue given in Proposition \ref{singleres}(a), as a function of $s_3$, is {\it analytic} at $s_3=a_m-\delta_3$, for $\delta_3\in\Z_{\ge0}$.   (Note that the factor $\gg{s_3-a_m}$ is ``missing'' from our above formula for this residue.)   
  
    \item  The residue given in Proposition \ref{singleres}(b), as a function of $s_3$, is {\it analytic} at $s_3=a_m-\delta_3$ and $s_3=a_n-\delta_3$, for $\delta_3\in\Z_{\ge0}$.   (Note that the factors  $\gg{s_3-a_m}$ and $\gg{s_3-a_n}$ are ``missing'' from our above formula for this residue.) Similarly, as a function of $s_1$, this residue is analytic at $s_1=-a_p-\delta_1$ and $s_1=-a_q-\delta_1$, for  $\delta_1\in\Z_{\ge0}$.\end{itemize}
    
    In light of the above observations, it will suffice to compute the following ``two-variable-at-a-time'' residues.     
\begin{prop}  \label{doubleres}
\noindent{\rm (a)} For   $ \delta_1,\delta_2\in\Z_{\ge0} $,  we have\begin{align}&\hbox{\rm Res}_{s_1=-a_1-\delta_1}\bigl(\res{s_2=-a_1-a_4-\delta_2}{\mhat{4}{a}{s}}\bigr)\nn
    \\& =\frac{(-1)^{\delta_1}  \gg{a_4-a_1-\delta_1}}{\delta_1!} \biggl[\prod_{ k=2}^3\gg{s_3-a_k}  \gg{a_k-a_1-\delta_2}  \gg{a_k-a_4-\delta_2}\biggr]    \nn     \\& \times p_{\delta_2}\biggl[{-\delta_1, a_2-a_4-\delta_2,s_3-a_3;\atop 1+a_1-a_3, a_2-a_1-\delta_1-\delta_2,s_3-a_4-\delta_2}\biggr]=  \gg{a_4-a_1-\delta_1}\nn
    \\&\times \biggl[\prod_{ k=2}^3\gg{s_3-a_k}  \gg{a_k-a_1-\delta_1}  \gg{a_k-a_4-\delta_2}\biggr]     f_{\delta_1,\delta_2}(s_3,a),\label{s1s2-res}\end{align}where $f_{\delta_1,\delta_2}$ is a polynomial of degree at most $2\delta_1+\delta_2$.

\noindent{\rm(b)} For  $ \delta_1,\delta_3\in\Z_{\ge0} $,  we have\begin{align}&\res{s_3=a_2-\delta_3}\bigl(\res{s_1=-a_1-\delta_1}{\mhat{4}{a}{s}})\nn=\frac{(-1)^{\delta_3}\gg{a_2-a_1-\delta_1}\gg{s_2+a_1+a_2}}{\delta_3! \gg{s_2 +a_1+a_2+\delta_1-\delta_3}}\\&\times\biggl[\prod_{k=3}^4\gg{a_k-a_1-\delta_1}\gg{a_2-a_k-\delta_3}
\gg{s_2+a_1+a_k}\biggr]\nn
\\&\times p_{\delta_1}\biggl[{-\delta_3,a_2-a_3-\delta_3,a_2-a_4-\delta_3;\atop s_2 +a_1+a_2-\delta_3, 1+a_1+a_2-s_2-\delta_3,a_2-a_1-\delta_1-\delta_3}\biggr]\nr&= \gg{a_2-a_1-\delta_1}  \biggl[\prod_{k=3}^4\gg{a_k-a_1-\delta_1}\gg{a_2-a_k-\delta_3}
\gg{s_2+a_1+a_k}\biggr]\nr&\times g_{\delta_1,\delta_3}(s_2,a),\end{align} where $g_{\delta_1,\delta_3}$ is a polynomial of degree at most $2\delta_1+\delta_3$.

\noindent{\rm(c)} For   $ \delta_2,\delta_3\in\Z_{\ge0} $,  we
   have\begin{align}&\res{s_3=a_3-\delta_3}\bigl(\res{s_2=-a_1-a_4-\delta_2}{\mhat{4}{a}{s}}\bigr)\nn
    \\& = \frac{(-1)^{\delta_3}\gg{a_3-a_2-\delta_3}}{\delta_3!}\biggl[\prod_{j\in\{1,4\} }\gg{s_1+a_j}   \gg{a_2-a_j-\delta_2} \gg{a_3-a_j-\delta_2}\biggr]    \nn     \\& \times p_{\delta_2}\biggl[{s_1+a_1, a_2-a_4-\delta_2,-\delta_3;\atop 1+a_1-a_3,s_1+a_2-\delta_2,a_3-a_4-\delta_2-\delta_3}\biggr] = \gg{a_3-a_2-\delta_3} 
      \nr&\times\biggl[\prod_{j\in\{1,4\} }\gg{s_1+a_j}   \gg{a_2-a_j-\delta_2} \gg{a_3-a_j-\delta_3}\biggr]  h_{\delta_2,\delta_3}(s_1,a), \end{align}where $h_{\delta_2,\delta_3}$ is a polynomial of degree at most $\delta_2+2\delta_3$.
    \end{prop}
\begin{pf} We prove part (a) only; proofs of the other parts are similar.  
The first equality in \eqref{s1s2-res} follows immediately from Proposition \ref{singleres}(b).   We rewrite this equality as follows:
\begin{align}&\hbox{\rm Res}_{s_1=-a_1-\delta_1}\bigl(\res{s_2=-a_1-a_4-\delta_2}{\mhat{4}{a}{s}}\bigr)\nn
    \\& = \gg{a_4-a_1-\delta_1}  \biggl[\prod_{ k=2}^3\gg{s_3-a_k}  \gg{a_k-a_1-\delta_1}  \gg{a_k-a_4-\delta_2}\biggr]\cdot\biggl\{\frac{(-1)^{\delta_1}  }{\delta_1!}    \nn     \\& \times\biggl[\prod_{ k=2}^3\frac{ \gg{a_k-a_1-\delta_2} }{\gg{a_k-a_1-\delta_1}}\biggr]   p_{\delta_2}\biggl[{-\delta_1, a_2-a_4-\delta_2,s_3-a_3;\atop 1+a_1-a_3, a_2-a_1-\delta_1-\delta_2,s_3-a_4-\delta_2}\biggr] \biggr\}.\label{s1s2-res2}\end{align} Denote the quantity in large braces, in \eqref{s1s2-res2}, by $f_{\delta_1,\delta_2}(s_3,a)$:  then we wish to show that $f_{\delta_1,\delta_2}(s_3,a)$ is a polynomial of degree at most $2\delta_1+\delta_2$.
    We consider two cases:  (i) $\delta_1>\delta_2$ and (ii) $\delta_1\le \delta_2$.  In the first case, we note that $$\prod_{ k=2}^3\frac{ \gg{a_k-a_1-\delta_2} }{\gg{a_k-a_1-\delta_1}}=\prod_{ k=2}^3(a_k-a_1-\delta_1)_{\delta_1-\delta_2}$$is a polynomial of degree $
    2(\delta_1-\delta_2)$.  Since $p_{\delta_2}$ has degree at most $3\delta_2$, it follows that $f_{\delta_1,\delta_2}(s_3,a)$ has degree at most $ 2(\delta_1-\delta_2)+3\delta_2=2\delta_1+\delta_2$, as claimed.
    On the other hand, if $ \delta_1\le\delta_2$, then  by part (c) of Lemma \ref{pdelta-lem}, there is a polynomial $q_{\delta_1,\delta_2}$ such that 
    \begin{align}&p_{\delta_2}\biggl[{-\delta_1, a_2-a_4-\delta_2,s_3-a_3;\atop 1+a_1-a_3, a_2-a_1-\delta_1-\delta_2,s_3-a_4-\delta_2}\biggr] \nr&=(1+a_1-a_3+\delta_1)_{\delta_2-\delta_1}(a_2-a_1 -\delta_2)_{\delta_2-\delta_1} \,q_{\delta_1,\delta_2}(s_3,a)
     \nr&=\frac{\gg{1+a_1-a_3+\delta_2}\gg{a_2-a_1 -\delta_1}}{\gg{1+a_1-a_3+\delta_1}\gg{a_2-a_1 -\delta_2}} \,q_{\delta_1,\delta_2}(s_3,a)
          \nr&=(-1)^{\delta_2-\delta_1}\frac{\gg{a_3-a_1-\delta_1}\gg{a_2-a_1 -\delta_1}}{\gg{a_3-a_1-\delta_2}\gg{a_2-a_1 -\delta_2}}\, q_{\delta_1,\delta_2}(s_3,a),\label{q}\end{align}The last step because $\gg{s}\gg{1-s}=\pi /\sin(\pi s)$.  The quantity in large braces, in \eqref{s1s2-res2}, then equals$$\frac{(-1)^{\delta_2 }}{\delta_1!} q_{\delta_1,\delta_2}(s_3,a).$$The degree of $q_{\delta_1,\delta_2}$ is, by \eqref{q}, less than or equal to $3\delta_2-2(\delta_2-\delta_1)=2\delta_1+\delta_2$.    So part (a) of our proposition is proved.\hfill  $\Box$\end{pf}  

\begin{rem}{\rm The above information on the maximum degrees of the  polynomials $f_{\delta_1,\delta_2}$, $g_{\delta_1,\delta_3}$, and $h_{\delta_2,\delta_3}$ will be critical to the calculations in \cite{GW}.}\end{rem}

Finally, we compute our three-variable residue.  This may be accomplished by computing the appropriate residue of any of the three parts of Proposition \ref{doubleres}.  We choose to begin with part (a) of this proposition. (Other approaches will yield results that look different, {\it a priori}, but can be transformed into each other  by various functional equations for $p_\delta$.)
\begin{prop}   Let  $ \delta_1,\delta_2,\delta_3\in\Z_{\ge0} $.    Then
\begin{align*}&\hbox{\rm Res}_{s_3=a_3-\delta_1}\bigl(\hbox{\rm Res}_{s_1=-a_1-\delta_1}\bigl(\res{s_2=-a_1-a_4-\delta_2}{\mhat{4}{a}{s}}\bigr)\bigr)=\frac{(-1)^{\delta_1+\delta_3} }{\delta_1!\delta_3!} 
\nn
    \\& \times \gg{a_4-a_1-\delta_1} \gg{a_3-a_2-\delta_3} 
  \biggl[\prod_{ k=2}^3  \gg{a_k-a_1-\delta_2}  \gg{a_k-a_4-\delta_2}\biggr]    \nn     \\& \times p_{\delta_2}\biggl[{-\delta_1, a_2-a_4-\delta_2,-\delta_3;\atop 1+a_1-a_3, a_2-a_1-\delta_1-\delta_2,a_3-a_4-\delta_2-\delta_3}\biggr].\end{align*}\end{prop}
 
\begin{pf} This follows immediately from Proposition \ref{doubleres}(a).\hfill  $\Box$\end{pf} 
\end{document}